\input amstex
\documentstyle{amsppt}
\NoBlackBoxes
\topmatter
\title
A note on solvable Lie groups without lattices and the F\'elix-Thomas models of fibrations
\endtitle
\author
Aleksy Tralle
\endauthor
\affil University of Warmia and Mazuria in Olsztyn
\endaffil
\address \D Zo\l\/nierska 14A, 10-561 Olsztyn, Poland
\endaddress
\email
tralle\@matman.uwm.edu.pl
\endemail
\abstract
In this article we show that  solvable Lie group $G^{bg}$ constructed in the paper of Benson and Gordon \cite{BG} has no lattices. This yields a negative answer to a question posed by several authors \cite{BG, FLS, TO} in the context of symplectic geometry. The result is of conceptual interest, because it is obtained by a delicate application of the F\'elix-Thomas theorem on models of non-nilpotent fibrations.
\endabstract
\endtopmatter
\rightheadtext{Solvable Lie groups without lattices}
\document
\head 1. Introduction
\endhead
\vskip6pt
Aspherical manifolds play an important role in many areas of geometry. For example, in symplectic geometry, such manifolds are used as a basic source of {\it symplectically aspherical} manifolds $(M,\omega)$, i.e. symplectic manifolds such that 
$$\int_{S^2}f^*\omega=0$$
for every map $f: S^2\to M$. Such manifolds have explicitly appeared in papers of Floer \cite{F} and Hofer \cite{H} in context of Lagrangian intersection theory. In these papers, the authors used an analytic advantage of symplectic asphericity: namely, it excludes the appearence of non-trivial pseudo-holomorphic spheres in $M$. Since symplectic asphericity can be equivalently described by the condition $\omega|_{\pi_2(M)}=0$, the role of asphericity becomes clear. For example, this condition was an important ingredient in proving of the Arnold conjecture about symplectic fixed points \cite{OR}, and, also, it was used in estimating the number of closed trajectories of symplectic magnetic fields in \cite{RT}. Note that explicit examples of symplectically aspherical manifolds with non-trivial $\pi_2(M)$  appeared very recently \cite{G}.  
\vskip6pt
Solvmanifolds, i.e. factor-spaces of solvable Lie groups $G$ over discrete co-compact subgroups (lattices) $\Gamma$ constitute a large and important subclass of the class of aspherical manifolds. This subclass is interesting not only on its own right, but very often serves as a source of testing examples in many areas of geometry (see, e.g. \cite{O}).
\vskip6pt
Formally, the present paper contains one explicit example of a completely solvable Lie group without lattices, answering (in a negative way) to a question posed in \cite{BG, FG, FLS, TO}. In greater detail,  
in \cite{BG} the following completely solvable Lie algebra was constructed:
$$\frak g^{bg}=\{S,T, X_1, Y_1, Z_1, X_2, Y_2, Z_2\}$$
$$[X_1,Y_1]=Z_1,\,\,[X_2,Y_2]=Z_2$$
$$[S,X_1]=X_1, [S,X_2]=-X_2,\,\,[S,Y_1]=-2Y_1,\,\,[S,Y_2]=2Y_2$$
$$[S,Z_1]=-Z_1,\,\,[S,Z_2]=Z_2.$$
Here $\frak g^{bg}$ is defined by the Lie brackets of the base elements, and all the brackets which are not written, are assumed to be zero.
\vskip6pt 
It was noted in \cite{BG} that this algebra is unimodular (i.e. $\operatorname{Tr}ad\,X=0$, for all $X\in\frak g^{bg}$).  Hence, the corresponding simply connected Lie group $G$ satisfies the necessary condition for the existence of a lattice. However, the authors pointed out that it is not known if there is a lattice in $G$ or not. If $G$ admitted a lattice, one would obtain an example of (aspherical) closed symplectic manifold which is formal but non-Lefschetz (this observation was made by Fern\'andez, de Le\'on and Saralegui \cite{FLS}, cf. also \cite{AFGM}).  It is not known if such {\it aspherical} closed symplectic manifolds exist. The importance of constructing them was discussed in \cite{IRTU}. In general, relations between formality and hard Lefschetz property are not yet clear (cf. \cite{AFGM, FLS, IRTU, M}). For example, there is a conjecture in \cite{M} that any closed symplectic manifold satisfying the hard Lefschetz property must be formal. Our result indicates the difficulties of constructing examples of that kind in the class of solvmanifolds.  
\vskip6pt
Although  we give an answer to a  question posed in \cite{BG, FG, TO}, the paper contains a result also of conceptual interest. There are several reasons for that. The first  reason is that it is (in general)  difficult to find whether the given solvable (non-nilpotent)  Lie group contains a lattice,even using  the known criteria \cite{A, S2, VGS}.  Hence, it is  important (in various geometric contexts) to have at least examples of non-nilpotent solvable Lie groups with and without lattices. 
\vskip6pt
Here is the second reason why the result is important. It is obtained by a delicate application of  rational homotopy theory. In greater detail, if the given Lie group $G^{bg}$ contained a lattice $\Gamma$, we could associate with $G/\Gamma$ the {\it Mostow fibration}
$$N/\Gamma\cap N \to G/\Gamma\to G/\Gamma N$$
and calculate the minimal model of the total space $G/\Gamma$ in two different ways. The first of these uses the Chevalley-Eilenberg complex $(\Lambda\frak g^*,\delta)$, and the second one -- the F\'elix-Thomas refinement of the method of KS-extensions \cite{FT}. Since these two methods yield different results, we obtain a contradiction. This neat application of the F\'elix-Thomas theorem is of independent interest (cf. \cite{OT, K}). Moreover, the methods developed in the paper seem to give a method of explicit constructing solvable Lie groups without lattices.
\vskip6pt
The paper is organized as follows. The main result is given as a Theorem in Section 3. In this Section, the proof of the Theorem is given. However, since there is an essential amount of direct but lengthy calculation, we postpone it to Section 4, in order to prevent the details of obscuring the idea. Also, we need an explicit construction of the bigraded model of a commutative graded algebra from \cite{HS}. We give a brief summary of this construction in Section 5.
\vskip6pt
\subhead Acknowledgement
\endsubhead
The author expresses his thanks to Marisa Fern\'andez, Yuli Rudyak and Sergei Merkulov for inspiring discussions and to John Oprea and Vladimir Gorbatsevich for reading the manuscript and valuable criticism and advice. This research was partially supported by the Polish Research Committee (KBN). 
\head 2. Preliminaries
\endhead
The result of this paper involves solvable Lie groups of a special type. These are called {\it completely solvable} Lie groups. We call a Lie group $G$ completely solvable, if the corresponding Lie algebra $\frak g$ satisfies the property that any adjoint linear operator $ad\,S: \frak g\to\frak g,\,\,S\in\frak g$ has only real eigenvalues. It can be shown that any completely solvable Lie group is solvable. Note that, obviously, any nilpotent Lie group is completely solvable.
\vskip6pt
To get the main theorem of this article, we use some techniques from rational homotopy theory. However, since now there are many good books and surveys on this topic, we restrict ourselves to references \cite{FH, H, L, Ta, TO, T}. We use freely the notion of the {\it minimal model} of a space $X$ and denote it $\Cal M(X)$. The minimal model of a differential graded algebra, say $(\Cal A,d)$ is denoted as $\Cal M(\Cal A)$. Recall that $\Cal M(X)$ (resp. $\Cal M(\Cal A)$ is a free minimal differential graded algebra, i.e. a tensor product of an exterior and a symmetric graded algebra, equipped with a differential. We use notation, reflecting this fact, by writing $\Cal M(X)=(\Lambda V, d_V)$, where $V$ is a graded vector space of free generators. Sometimes, we use the term "model" of a space $X$, and, in this case, we understand this simply as a differential graded algebra, say $\Cal A(X)$ such that it admits a morphism $\Cal A(X)\to A(X)$ inducing isomorphism on the cohomology level. In the sequel, we denote by $A(X)$ the differential graded algebra of polynomial differential forms. Note that all differential graded algebras are considered over $\Bbb Q$.
\vskip6pt
Also, there is a way to model fibrations (the algebraic models of the latter are called the {\it Koszul-Sullivan models} (or, briefly, K.S. models). In greater detail, let $F\to E\to B$ be any fibration of path connected spaces. Then there exists a commutative diagram of differential graded algebras
$$
\CD
(A(B),d_B) @>>> (A(E),d_E) @>>> (A(F),d_F)\\
@A{\psi_B}AA @A{\psi_E}AA @A{\psi_F}AA\\
(\Lambda\,V,d_V) @>>> (\Lambda(V\oplus Y),D) @>>>(\Lambda\,Y,d_Y)
\endCD
$$
where $\psi_B$ and $\psi_E$ induce cohomology isomorphisms and there is a well ordered index set $J$ such that
$$D(y_{\alpha})\in \Lambda\,V\otimes\Lambda\,Y_{<\alpha},\,\,\text{and}\,\,\,\alpha<\beta\implies |y_{\alpha}|<|y_{\beta}|.$$
Here, $y_{\alpha},\alpha\in J$ constitute a basis for $Y$ and notation $\Lambda\,Y_{<\alpha}$ is used for a subalgebra in $\Lambda\,Y$ generated by all $y_{\beta}\in Y, \beta<\alpha$.
\vskip6pt
Note that $\psi_F$ does not necessarily induce a cohomology isomorphism, so that $(\Lambda\,Y,d_Y)$ is not necessarily a minimal model of $F$. In the case when $\pi_1(B)$ acts nilpotently on the cohomology of the fiber, the well-known theorem of Grivel-Halperin-Thomas shows that $(\Lambda\,Y,d_Y)$ is, in fact, a minimal model of $F$. If the action of $\pi_1(B)$ is not nilpotent (and this is the situation in the present paper), we need the following refinement of the cited result.
\proclaim{F\'elix-Thomas Theorem} \cite{FT} Let $F\to E\to B$ be a fibration and let $U$ denote the largest $\pi_1(B)$-submodule of $H^*(F)$ on which $\pi_1(B)$ acts nilpotently. Assume that the following conditions are satisfied:
\roster
\item "(i)" $H^*(F)$ is a vector space of finite type;
\item "(ii)" $B$ is a nilpotent space (i.e. $\pi_1(B)$ is a nilpotent group which acts nilpotently on $\pi_i(B)$ for all $i\geq 1$).
\endroster
Then, in the K.S. model of the fibration
$$
\CD
A(B) @>>> A(E) @>>> A(F)\\
@AAA @AAA @A{\psi}AA\\
(\Lambda V,d_V) @>>> (\Lambda(V\oplus Y),D) @>>> (\Lambda Y,d_Y)
\endCD
$$
the differential graded algebra homomorphism $\psi: (\Lambda Y,d_Y)\to A(F)$ induces an isomorphism
$$\psi^*: H^*(\Lambda Y,d_Y)\cong U.$$
\endproclaim
Note that for a nilpotent fibrations, the original Grivel-Halperin-Thomas theorem is recovered.
\vskip6pt
In the sequel, we will need the fact, that solvmanifolds $G/\Gamma$ of completely solvable Lie groups $G$ have a particular free model, given by the Nomizu--Hattori theorem \cite{TO}. This theorem shows that the standard Chevalley--Eilenberg complex for the Lie algebra cohomology $H^*(\frak g)$ is also a free model for the manifold $G/\Gamma$. Recall that the Chevalley--Eilenberg complex for the Lie algebra cohomology is defined as follows. One considers the differential graded algebra of the form
$$(\Lambda\frak g^*,\delta)=(\Lambda(x_1,..., x_n),\delta)$$
with $|x_i|=1, i=1,...,n$ and the differential $\delta$  given by the formula
$$\delta x_i(e_k,e_r)=x_i([e_k,e_r]),$$
where $x_i$ denote the vectors dual to the vectors of the base, say, $\{e_1,...,e_n\}$ of the Lie algebra $\frak g$. Of course, this differential graded algebra calculates the cohomology of $\frak g$ as well as the cohomology of the manifold $G/\Gamma$, but we have also a more delicate result:
$$\Cal M(G/\Gamma)\cong \Cal M(\Lambda\frak g^*,\delta).$$
Hence, we can calculate the minimal model of a completely solvable solvmanifold from the Lie algebra data in a purely algebraic way. This will be used in the sequel.
\vskip6pt
Finally, in the proof of the main result of this paper, we need a special way of constructing free models of algebras and spaces, given in the paper of Halperin and Stasheff \cite{HS}. We refer to this paper for details. Here we only mention that we use this method explicitly, and, therefore, we outline it in a separate Section 5 for the convenience of the reader. Here we only mention that to get such free model, we need the graded commutative algebra $H$ (in case of spaces, $H$ is a cohomology algebra). We endow $H$ with zero differential and construct in a canonical way a {\it bigraded model} $(\Lambda Z,d)$ from the given $(H,0)$. Then, there is a free model (not minimal in general) $(\Lambda Z,D)$ which is the same as a graded algebra but with a new ("perturbed") differential $D$:
$$D=d_0+d_1+ d_2+..., \,\,\, D-d: Z_n\to F_{n-2}(\Lambda\,Z),$$
where $F_n(\Lambda\,Z)=\sum_{m\leq n}(\Lambda\,Z)_m$ denote the spaces of a canonical filtration on $\Lambda\,Z$. 
\vskip6pt
In Section 3 we consider the 3-dimensional Heisenberg Lie algebra, i.e. the Lie algebra $\frak n_3=\{e_1,e_2,e_3\}$ with Lie brackets given by $[e_1,e_2]=e_3, \,[e_i,e_j]=0$ for all other pairs $e_i, e_j$.
  
\head 3.  Result \endhead

Let us formulate the main result of this work.
\proclaim{Theorem} The solvable Lie group $G^{bg}$ has no lattices.
\endproclaim
The rest of this Section is devoted to the proof of this claim. However, we need a preparatory work, which we summarize in the forthcoming Lemmas 1-5.
\vskip6pt 
Let $G$ be a solvable Lie group. Denote by $N$ its nilradical, i.e. a maximal connected nilpotent Lie subgroup.
\proclaim{Lemma 1} Any simply connected completely solvable Lie group $G$ containing a lattice and satisfying the condition $\dim\,G/N=1$ is splittable into a semidirect product of the form $G=\Bbb R\rtimes_{\varphi}N$.
\endproclaim
\demo{Proof} Let $\Gamma$ be a lattice in $G$. Put $\Gamma_N=\Gamma\cap N$. It is well known that $\Gamma_N$ is a lattice in $N$ \cite{R}. In general, there is the following commutative diagram with exact horizontal rows
$$
\CD
1 @>>> N @>>> G @>{\tilde\pi}>> G/N @>>> 1\\
@AAA @AAA @AAA @AAA @AAA\\
1 @>>> \Gamma_N @>>> \Gamma @>{\pi}>> \Gamma/\Gamma_N @>>> 1
\endCD
$$
such that $G/N\cong \Bbb R^s$ and $\Gamma/\Gamma_N\cong \Bbb Z^s$.
Note that in our case $s=1$. Since $\Gamma/\Gamma_N\cong \Bbb Z$ is free, the lower row in the diagram admits a splitting $s: \Bbb Z\to \Gamma$.  It follows that $\Gamma\cong\Bbb Z\rtimes_{\psi}\Gamma_N$. 
\vskip6pt
Note that $G/N$ and $G$ are completely solvable Lie groups, and one can apply to them the following result of Sait\^o \cite{S1}:
\proclaim{Sait\^o's Theorem} For any simply connected completely solvable Lie groups $G$ and $G'$ such that $G$ contains a lattice $\Gamma$, any homomorphism $\alpha: \Gamma\to G'$ can be uniquely extended to a homomorphism $\bar\alpha: G\to G'$.
\endproclaim
Hence, we extend the given section $s: \Bbb Z\cong \Gamma/\Gamma_N\to \Gamma\subset G$ to the map $\tilde s: G/N\to G$. We claim that $\tilde s$ is a section of the first row in the commutative diagram above. Indeed, $\tilde\pi\circ\tilde s: G/N\to G/N$ is an extension of the homomorphism $\pi\circ s: \Gamma/\Gamma_N\to \Gamma$:
$$\tilde\pi\circ\tilde s|_{\Gamma/\Gamma_N}=\tilde\pi|_{s(\Gamma/\Gamma_N)}=\pi\circ s=\text{id}_{\Gamma/\Gamma_N}$$
Applying the Sait\^o theorem again (i.e. the uniqueness of the extension), we get $\tilde\pi\circ\tilde s=\text{id}_{G/N}$, as required. It follows that $G$ is splittable, and the proof is completed. 
\enddemo
\remark{Remark} Note that from the construction $\varphi$ restricted to ${\Bbb Z}$ equals $\psi.$
Also, the proof cannot be extended to groups with $\dim\,G/N>1$, since the group $\Bbb Z^k, k>1$ is not free.
\endremark 
\vskip7pt
Now, if $G=A\rtimes_{\theta}N$ is completely solvable with $\dim\,A=1$, we have (in general) two different splittings:
$$G=A\rtimes_{\theta}N \quad \text{and} \quad G=A_1\rtimes_{\varphi}N$$
and the second one comes from the Lemma. The splitting $G=A_1\rtimes_{\varphi}N$ satisfies the property $\varphi(\Bbb Z^{A_1})(\Gamma_N)\subset \Gamma_N$. Hence it is natural to call the first splitting {\it initial}, while the second one -- {\it compatible} (with the splitting of lattice $\Gamma$). Note that of course $A\cong A_1\cong \Bbb R$ and $\Bbb Z^{A_1}\cong \Bbb Z^A\cong \Bbb Z$ naturally denote the corresponding lattices.
\vskip6pt
Consider the Mostow fibration corresponding to a solvmanifold $G/\Gamma$ such that $G$ satisfies the conditions of the Lemma:
$$N/\Gamma_N\to G/\Gamma\to G/N\Gamma\cong S^1$$
In the sequel we will need a description of the $\pi_1(S^1)\cong\Bbb Z$-action on the cohomology of the fiber $H^*(N/\Gamma_N)\cong H^*(\frak n)$ determined by this fibration. 
\proclaim{Lemma 2} Let $G/\Gamma$ be a compact completely solvable solvmanifold such that $\dim\,G/N=1$. Then, for the corresponding Mostow fibration $\pi_1(S^1)$-action on $H^*(N/\Gamma_N)$ is given by:
\roster
\item taking the compatible splitting $A_1\rtimes_{\varphi}N$,
\item restricting $\bar\varphi=d\varphi: A_1\to \operatorname{Aut}\,(\frak n)$ to $\psi: \Bbb Z^{A_1}\to \operatorname{Aut}\,(\frak n)$,
\item taking the dual automorphism $\psi^t:\Bbb Z^{A_1}\to \operatorname{Aut}(\frak n^*)$,
\item extending it to an automorphism of the exterior algebra $\Lambda\psi^t:\Bbb Z^{A_1}\to \operatorname{Aut}(\Lambda\frak n^*)$ as differential graded algebras,
\item and taking the induced automorphisms on cohomology $(\Lambda\psi^t)^*:\Bbb Z^{A_1}\to \operatorname{Aut}(H^*(\Lambda\frak n^*,\delta))$.
\endroster
\endproclaim
\demo{Proof} From Lemma 1 and the compatibility of the splitting we conclude that $G/\Gamma$ satisfies the conditions of Theorem 2.4 in \cite{OT}. Applying this theorem we complete the proof.
\enddemo
\vskip6pt
\remark{Remark} In \cite{OT} $\bar\varphi=d\varphi$ is denoted simply by $\varphi$ since it does not cause any confusion. In this paper it is more convenient to distinguish $\varphi: A_1\to \operatorname{Aut}(N)$ and $\bar\varphi: A_1\to \operatorname{Aut}(\frak n)$.
\endremark
\vskip6pt 
In the proof of the next Lemma we will need two equivalent interpretations of a semidirect product of two Lie groups, say, $G_1$ and $G_2$ (see \cite{VO}). In the first interpretation, $G=G_2\rtimes_b G_1$, and we understand $G$ as a Cartesian product $G_2\times G_1$ with the multiplication determined by some homomorphism $b: G_2\to \operatorname{Aut}(G_1)$ according to the formula
$$(g_2,g_1)\cdot (h_2,h_1)=(g_2h_2,g_1b(g_2)h_1).$$
According to the second interpretation, $G=G_2\cdot G_1$ (as a product of subgroups $G_1$ and $G_2$, $G_1\cap G_2=\{e\}$, with $G_1$ being normal). In that case, we get a homomorphism $g: G_2\to\operatorname{Aut}(G_1)$ determined by the expression:
$$b(g_2)g_1=g_2g_1g_2^{-1}, \,\,g_1\in G_1, \,g_2\in G_2.$$
\proclaim{Lemma 3} Assume that $G$ is completely solvable, admits a lattice $\Gamma$ and has a splitting $G=A\rtimes_{\theta}N$, such that $\dim\,A=1$. Consider the compatible splitting $G=A_1\rtimes_{\varphi}N$ given by Lemma 1 and the Mostow fibration
$$N/\Gamma\cap N \to G/\Gamma \to G/N\Gamma.$$
Then the $\pi_1(S^1)$-action on $H^*(N/\Gamma_N)$ is expressed through initial homomorphism $\theta$ as follows:
$$(\Lambda\psi^t)^*(a_1)=(\Lambda\bar\theta^t)^*(a)\circ (\Lambda\alpha^t)$$
where the homomorphism $\psi$ comes from the compatible splitting, and $\alpha\in\operatorname{Aut}(\frak n)$ is a unipotent automorphism of the Lie algebra $\frak n$. Here $a_1\in \Bbb Z^{A_1}\subset A_1$ is a generator of $\Bbb Z^{A_1}\cong\Bbb Z$ and $a\in A$ is an element uniquely determined by $a_1$.
\endproclaim
\demo{Proof} Consider the given splittings
$$G=A_1\rtimes_{\varphi}N,\quad \text{and}\quad G=A\rtimes_{\theta}N.$$
Note that we think about both of them as decompositions 
$$G=A_1\cdot N,\quad G=A\cdot N, \quad A_1\cap N=A\cap N=\{1\}.$$
of the same group $G$ into two different products of {\it subgroups}, with the same normal term $N$ (see the remarks before Lemma 3). Using this interpretation of the semidirect product, we can assume that
$$\varphi(a_1)(n)=a_1\cdot n\cdot a_1^{-1},\,\,\, \theta(a)=a\cdot n\cdot a^{-1},\,\,\,\text{for any}\,\,\, a\in A,\,a_1\in A_1.$$
Now, take $a_1\in\Bbb Z^{A_1}$. Since $a_1\in A_1\subset G$, we can decompose it according to the initial spliting:
$$a_1=a\cdot n_{a_1},\,\, a\in A,\,\, n_{a_1}\in N.$$
Applying the formulas above we get
$$\varphi(a_1)(n)=a_1\cdot n\cdot a_1^{-1}=(a\cdot n_{a_1})\cdot n\cdot (a\cdot n_{a_1})^{-1}=a( n_{a_1}\cdot n\cdot n_{a_1}^{-1})a^{-1}=\theta(a)\circ I(n_{a_1})(n),$$ 
where $I(n_{a_1})$ denotes the inner automorphism of $N$ determined by $n_{a_1}$. Finally
$$\varphi(a_1)=\theta(a)\circ I(n_{a_1}).$$
Note that obviously $\alpha=Ad(n_{a_1})=I(n_{a_1})_*$ is a unipotent automorphism of $\frak n$, and this remark completes the proof.
\enddemo
\remark{Remark} The fact that $Ad(n_{a_1})$ is unipotent, follows from the well-known formula $Ad(exp\,X)=e^{ad\,X}$ for all $X$ from the Lie algebra of any Lie group, the fact that $\exp$ is a diffeomorphism for nilpotent Lie groups and the nilpotency of any $ad\,X, X\in\frak n$ (the latter is Engel's theorem).
\endremark
\vskip7pt
Let $\frak g^{bg}$ be the Benson-Gordon Lie algebra and let $G^{bg}$ denote the simply-connected completely solvable Lie group whose Lie algebra is $\frak g^{bg}$.  
\proclaim{Lemma 4} For the Lie algebra $\frak g^{bg}$ the following holds:
\roster
\item $\frak g^{bg}$ is a semidirect sum
$$\frak g^{bg}=\langle S\rangle\oplus_{\beta}\langle T, X_1,Y_1,Z_1,X_2,Y_2,Z_2\rangle=\frak a\oplus\frak n$$
where $\frak a=\langle A\rangle$ is a one-dimensional subalgebra and $\frak n$ is a 7-dimensional nil-radical. $T$ belongs to the center of $\frak g$.
\item $\beta(tS)=ad\,tS$ for any $t\in\Bbb R$.
\item $$G^{bg}=\Bbb R\rtimes_{\theta}N, d\theta(t)=\bar\theta(t): \Bbb R\to \operatorname{Aut}(\frak n)$$
and, moreover, $\exp\beta(t)=\bar\theta(t)$, for any $t\in\Bbb R$
\endroster
In particular, $G^{bg}$ satisfies the assumptions of Lemmas 1-3.
\endproclaim
\demo{Proof} Obviously, $\frak g^{bg}$ is a semidirect sum of Lie algebras $\frak a$ and $\frak n$. The rest of the proof is a  description of the correspondence between semidirect sums of Lie algebras and simply connected Lie groups which are semidirect products (cf. \cite{OT, VO}).
\enddemo
\proclaim{Lemma 5} Let $G^{bg}=A\rtimes_{\theta}N, A\cong \Bbb R$ be the splitting given by Lemma 4. Let $\alpha\in \operatorname{Aut}(\frak n)$ be any automorphism of the form $\alpha=Ad(n), n\in N$. Then, for any $a\in A\cong \Bbb R$,  automorphism $d\theta(a)\circ\alpha: \frak n\to\frak n$ has triangular matrix in the base $\{A,B,X_1,Y_1,Z_1,X_2,Y_2,Z_2\}.$
\endproclaim
\demo{Proof} The proof is straightforward.  From Lemma 4, $$d\theta(a)=\bar\theta(a)=\exp\beta(a)=ae^{ad\,S}$$
Since
$$\beta=ad\,S=\text{diag}\,(0,1,-2,-1,-1,2,1),$$
we get, by exponentiating, $$\bar\theta(a)=\text{diag}(1,\nu,\nu^{-2},\nu^{-1},\nu^{-1},\nu^{2},\nu), \,\,\nu=e^a,a\in\Bbb R.$$

Now, consider $\bar\theta(a)\circ\alpha=\bar\theta\circ e^{ad\,R}: \frak n\to\frak n$, where $R$ is a linear combination of the vectors $T,X_1,Y_1,Z_1,X_2,Y_2,Z_2$.  
One can check that in {\it this} base $ad\,R: \frak n\to\frak n$ is triangular, with zeros on the diagonal. We have also proved that in the same base, $ae^{ad\,S}$ is diagonal as well, and the proof follows.
\enddemo 
Now we are ready to prove our result.
\subhead Proof of the Theorem
\endsubhead
Assume that $G^{bg}$ has a lattice. Then we can form the Mostow fibration
$$N/\Gamma_N\to G^{bg}/\Gamma\to G^{bg}/N\Gamma\cong S^1.$$
Moreover, from Lemma 4, $G/\Gamma$ satisfies the assumptions of Lemmas 2,3 and 5. By Lemma 3, the $\pi_1(S^1)$-action on $H^*(\Lambda\frak n^*,\delta)$ is given by $(\Lambda\bar\theta(a)^t)^*\circ(\Lambda\alpha^t)^*$. From the F\'elix-Thomas theorem, the model of the Mostow fibration is given by the sequence of the differential graded algebras
$$(\Lambda(a),0)\to (\Lambda(a)\otimes\Lambda\,V,D)\to (\Lambda\,V,d)$$
such that $H^*(\Lambda\,V,d)\cong U\subset H^*(\Lambda\frak n^*,\delta)$, where $U$ is a maximal nilpotent submodule with respect to $(\Lambda\bar\theta(a)^t)^*\circ (\Lambda\alpha^t)^*$-action. By Lemma 5, $\bar\theta(a)\circ\alpha$ has a triangular matrix. Let us write this condition explicitly, after the dualizing. Denote the dual base in $\frak n^*$ by $b,x_1,y_1,z_1,x_2,y_2,z_2$ and write the expressions
$$
\aligned
\eta(b)&=b\\
\eta(x_1)&=\nu x_1+\text{lin.comb.}\{y_1,z_1,x_2,y_2,z_2\}\\
\eta(y_1)&={1\over \nu^2}y_1+\text{lin.comb.}\{z_1,x_2,y_2,z_2\}\\
\eta(z_1)&={1\over \nu}z_1+\text{lin.comb.}\{x_2,y_2,z_2\}\\
\eta(x_2)&={1\over \nu}x_2+\text{lin.comb.}\{y_2,z_2\} \\
\eta(y_2)&=\nu^2y_2+\mu z_2\\
\eta(z_2)&=\nu z_2,\endaligned
$$
where $\eta$ stands for $\bar\theta^t\circ\alpha^t.$ Now we are calculating the induced $\Lambda\eta^*$-action on the cohomology $H^*(\Lambda\frak n^*)$. Note that, obviously, 
$$\frak n=\langle T\rangle\oplus\frak n_3\oplus\frak n_3,\,\,\text{and}\,\,H^*(\frak n)=\Lambda(b)\otimes H^*(\frak n_3)\otimes H^*(\frak n_3)$$
Since $\eta(T^*)=\eta(b)=b$, to calculate $U$, we need, in fact, only to analyze $H^*(\frak n_3)\otimes H^*(\frak n_3)$. First, we write the linear basis for the latter cohomology:
\vskip6pt
\roster
\item "$H^1$:" $[x_1],[y_1],[x_2], [y_2]$
\vskip6pt
\item "$H^2:$" $[x_1z_1],[y_1z_1], [x_1][x_2],,[x_1][y_2],[y_1][x_2]$,
\newline
 $[y_1][y_2],[x_2z_2],[y_2z_2]$
\vskip6pt
\item "$H^3$:" $[x_1z_1][y_1],[x_1z_1][x_2],[x_1z_1][y_2],[y_1z_1][x_2],[x_2z_2][x_1],$ 
\newline 
$[x_2z_2][y_1],[x_2z_2][y_2],[y_2z_2][x_1],[y_2z_2][y_1],[y_1z_1][y_2]$
\vskip6pt
\item "$H^4$:" $[x_1y_1z_1][x_2],[x_1y_1z_1][y_2],[x_2y_2z_2][x_1],[x_2y_2z_2][y_1],
[x_1z_1][x_2z_2],$ 
\newline
$[x_1z_1][y_2z_2],[y_1z_1][x_2z_2],[y_1z_1][y_2z_2]$
\vskip6pt
\item "$H^5$:" $[x_1y_1z_1][x_2z_2],[x_1y_1z_1][y_2z_2],[x_2y_2z_2][x_1z_1],[x_2y_2z_2][y_1z_1]$
\vskip6pt
\item "$H^6$:"  $[x_1y_1z_1][x_2y_2z_2]$
\endroster
\vskip6pt
The above formulas follow, since the cohomology $H^*(\Lambda\frak n_3)$ is obviously given by
$$
\aligned H^1(\Lambda\frak n_3^*)&=\langle [x_1],[y_1]\rangle\\
H^2(\Lambda\frak n_3^*)&=\langle [x_1z_1],[y_1z_1]\rangle\\
H^3(\Lambda\frak n_3^*)&=\langle [x_1y_1z_1]\rangle,\endaligned
$$
(and the same  holds for variables $x_2,y_2,z_2$), one easily obtains the above formulae for the cohomology $H^*(\Lambda\frak n^*)$.
\vskip6pt  
Once can check by direct calculation that in this base of $H^*$, the following holds:
\roster
\item "(i)" $\Lambda\eta^*$ has the triangular matrix,
\item "(ii)" the following vectors 
$$[x_1z_1],[x_1x_2],[y_1y_2],[x_2z_2],[y_1z_1][y_2z_2],[x_1y_1z_1][y_2], [y_1][x_2y_2z_2]\eqno (*)$$
satisfy the equality of the form
$$\Lambda\eta^*(w_j)=w_j+\text{lin. comb.}\{w_k, k>j\}.\eqno (**)$$
and these are the only base vectors with this property.
\endroster
Here, of course, $w_j$ stand for vectors $(*)$.
Although the calculation is straightforward, we give it in a separate section for the convenience of the reader. From (i) and (ii) we see that the maximal nilpotent submodule $U\subset H^*(\Lambda\frak n^*)$ is generated by elements $(*)$ and $b$. Indeed, if the matrix of a linear transformation, say $P: W\to W$ is triangular in some base $w_1,...,w_n$, then the subspace $W'$ satisfying the condition $W'=\{w\in W| (A-E)^m(w)=0,\,\text{for some}\,m\in\Bbb Z\}$, is generated by vectors of the form $(**)$.
\vskip6pt
It follows that 
$$U=\langle b,u_1,u_2,u_3,u_4,v_1,v_2,v_3\rangle$$
where $u_1=[x_1z_1],u_2=[x_2z_2],u_3=[x_1x_2],u_4=[y_1y_2], v_1=[y_1z_1][y_2][z_2],v_2=[x_1y_1z_1][y_2],v_3=[y_1][x_2y_2z_2]$. It follows that $U$ is a graded commutative algebra generated by elements $u_i, i=1,...,4$ of degree 2 and elements $v_j, j=1,...,3$ of degree 4 and one element $b$ of degree 1. The generators of even degrees satisfy the following relations:
$$\aligned 
u_i^2&=u_1u_3=u_2u_3=u_3u_4=0,\, i=1,...,4\\
u_1v_1&=u_1v_2=u_2v_1=u_3v_2=\\
u_3v_3&=u_4v_1=u_4v_2=u_4v_3=0\\
v_sv_t&=0\,\,\, \text{for all}\,\,\, s,t.\endaligned
$$
\vskip6pt
Finally, the F\'elix-Thomas theorem yields the following model of $G/\Gamma$:
$$\Cal A(G/\Gamma)\cong (\Lambda(a)\otimes\Lambda V,D)$$  
$$Da=0, H^*(\Lambda\,V,d)=U=\langle b, u_1,u_2,u_3,u_4,v_1,v_2,v_3\rangle$$
$$|a|=|b|=1,|u_i|=2,\,i=1,...,4, \,|v_1|=|v_2|=|v_3|=4.$$
for the {\it hypothetical} solvmanifold $G/\Gamma$.
\vskip6pt
Note that we can always assume that $(\Lambda\,V,d)$ is a minimal differential graded algebra obtained as follows: one calculates the minimal bigraded model $(\Lambda\,Z^{U},d_U)$ of $U$, takes some perturbation of the differential and obtains the free algebra $(\Lambda\,Z^U,D_U)$. Then, $(\Lambda\,V,d)$ is calculated as a minimal model of $(\Lambda\,Z^U,D_U)$:
$$(\Lambda\,V,d)\cong \Cal M(\Lambda\,Z^U, D_U).$$
We see that
$$(\Lambda\,Z^U,D_U)=(\Lambda(b),0)\otimes (\Lambda\,(Z^U_0\oplus Z^U_1\oplus...),D_U)$$
where $(Z^U_0)^2=\langle u_1,u_2,u_3,u_4\rangle$, $(Z^U_0)^4=\langle v_1,v_2,v_3\rangle$.
Further,
$$Z^U_1=(Z^U_1)^3\oplus (Z^U_1)^5\oplus (Z^U_1)^7$$
where
$$(Z^U_1)^3=\langle t_1,...,t_4\rangle,\,\,(Z^U_1)^5=\langle s_{ij}\rangle,\,\,(Z^U_1)^7=\langle q_{st}\rangle.$$
The differential $D_U$ on $Z^U_0\oplus Z^U_1$ is given as follows:
$$
\aligned
D_U(b)&=0\\ 
D_U(u_i)&=D_U(v_j)=0,\text{for all}\,\,i=1,...,4,j=1,...,3,\\
 D_U(t_i)&=u_i^2,\\
D_U(s_{ij})&=u_iv_j,\\
 D_U(q_{st})&=v_sv_t.
\endaligned
$$
The above expressions for the differential follow from the method of constructing of the bigraded model $(\Lambda\,Z^U,d_U)$ and the fact that $D_U-d_U: Z_n\to F_{n-2}(\Lambda\,Z^U)$, i.e. that the "perturbation" of the differential does not change it on the spaces $Z^U_0$ and $Z^U_1$.
 Finally,
$$\Cal A(G/\Gamma)=(\Lambda(a,b)\otimes\Lambda\,\bar V,D).$$
where $(\Lambda\,\bar V,d)\cong\Cal M(\Lambda\,(Z^U_0\oplus Z^U_1\oplus Z^U_2\oplus ...\oplus...), D_U)$ is a minimal differential graded algebra obtained by the procedure described above, and $Da=Db=0$ for the degree reasons.
For the same reason, the result is again minimal, and therefore
$$\Cal M(G/\Gamma)=(\Lambda(a,b)\otimes\Lambda\,\bar V,D), Da=Db=0.$$
Now, we want to get some conditions on the generators of $(\Lambda\,\bar V,d)$ coming from the fact that this algebra is obtained from the free algebra (although, in general, not minimal) $(\Lambda\,(Z^U_0\oplus Z^U_1\oplus Z^U_2\oplus...),D_U)$. To do this, we need to combine the method of constructing of $(\Lambda\,Z^U,D_U)$ as a perturbation of the bigraded model $(\Lambda\,Z^U,d_U)$ \cite{HS} and a method of calculating the minimal model of any free {\it simply-connected} and connected differential algebra given in \cite{L, Section 8}, or in \cite{TO,p. 158}. Note that we describe this methods separately, in Section 5, for the convenience of the reader. Here we omit it, in order to prevent the details of the construction from obscuring the idea of the proof.
\vskip6pt
The method of constructing the bigraded model $(\Lambda\,Z^U,d_U)$ shows that $(Z^U_i)^j\not=0,i>1$ only for $j>3$, since each space $(Z^U_i)$ is added to kill "new" unnecessary cocycles, which appeared at the previous step of the construction (see Section 5). For example,
$$(Z_2^U)^3=(H_1(\Lambda\,(Z_0\oplus Z_1), d_U)/H_1(\Lambda\,(Z_0\oplus Z_1)\cdot H^+_0(\Lambda\,(Z_0\oplus Z_1))^4$$
and we see that any new cocycle  must involve elements $t_i, s_{ij}$, or $q_{st}$, which necessarily raises the degree to at least 5, a contradiction. Now, as we have already mentioned, the non-minimality of $(\Lambda\,Z^U,D_U)$ can appear after the perturbation of the differential $d_U$, namely, 
$$D_U=(d_U)_0+(d_U)_1+\text{perturbed terms},$$
and in these perturbed terms, some linear expressions may appear.
\vskip6pt
 Applying now the method of calculating of $\Cal M(\Lambda\,Z^U,D_U)$, one can notice the following. First, in the perturbed model, {\it no expressions of the form
$$D_U(z_j)=v_j+\text{decomposables}$$  
can appear.} Indeed, since $|z_j|=3$, the previous remarks show that, otherwise, $z_j\in (Z^U_1)^3$, a contradiction. But since $v_j$ cannot be a linear term in any coboundary, the method of calculating the minimal model (see Section 5) yields
$$\Cal M(\Lambda\,(Z^U_0\oplus Z^U_1\oplus Z^U_2\oplus...),D_U)=(\Lambda\,\bar V,d),\,\, \bar V=\oplus_{ij}\bar V^j_i,\,\,V_0^4=\langle v_1,v_2,v_3\rangle,$$
since $\bar V_0^4\subset V'$ in the notation of the Proposition in Section 5.
\vskip6pt
 On the other hand, if $G$ contained $\Gamma$, the corresponding solvmanifold would be completely solvable, and as a consequence, it would have the same minimal model as $(\Lambda\frak g^*,\delta)$. The latter minimal model was calculated by Fern\'andez, de Le\'on and Saralegui \cite{FLS}. It was shown in this article that
$$\Cal M(\Lambda\frak g^*,\delta)\cong (\Lambda (a,b)\otimes\Lambda\,Z,d)$$
where $Z=Z_0\oplus Z_1\oplus Z_2$, $Z_0^2=\{b_i\}, Z_0^4=\{e\}, Z_1^3=\{c_i,g_i\}, Z_1^4=\{h_j\},Z_1^5=\{f\}, Z_2^4=\{x_i,y_j,z_j\}$.
Compairing the degrees of generators we see that 
$$\Cal M(G/\Gamma)\not\cong \Cal M(\Lambda\frak g^*,\delta),$$
since $\dim\,V_0^4=3\not=\dim\,Z^4_0=1$.
 Hence we have obtained a contradiction caused by the assumption that $G^{ab}$ admitted a lattice.

\head 4. Calculation of $U$
\endhead
Now, to simplify notation, we write  $\eta=\Lambda\bar\theta(a)^t\circ\Lambda\alpha^t$. Let us show, first, that $\eta^*: H^*(\Lambda\frak n^*)\to H^*(\Lambda\frak n^*)$ has a triangular matrix. Since $\eta^*$ is a morphism of degree $0$, it is sufficient to calculate the matrix of $\eta^*$ on each $H^i$ separately.
\subhead Case $H^1$
\endsubhead
Obviously, $\eta^*([x_1])=\nu [x_1]+\text{lin.comb.}\{[y_1],[x_2],[y_2]\}$.

Now,
$$\eta^*([y_1])=[\eta(y_1)]=[{1\over\nu^2}y_1+\text{lin. comb.}\{z_1,x_2,y_2,z_2\}].$$
We claim that $[\text{lin. comb.}\{z_1,x_2,y_2,z_2\}]=\text{lin. comb.}\{[x_2],[y_2]\}$. Indeed, one can write $[\mu_1z_1+\mu_2z_2+\delta_2x_2+\gamma_2y_2]=\text{lin. comb.}\{[x_2], [y_2]\}+[\mu_1z_1+\mu_2z_2]$. However, the expression $\mu_1z_1+\mu_2z_2$ can be a cocycle only if $\mu_1=\mu_2=0$.
\vskip6pt
In the same way we check the necessary conditions for the values of $\eta^*$ on the cohomology classes $[x_2], [y_2]$, and the matrix of $\eta$ restricted onto $H^1$ is triangular.
\subhead Case $H^2$
\endsubhead
This is done by the same straightforward calculation of $\eta^*[u]$ as $[\eta(u)]$ for any vector of the base of $H^2$. Note that to get the the triangular matrix, the vectors should be written in the same order as in the previous section.
\vskip6pt
For example, 
$$\eta^*([y_1x_2])=[({1\over\nu^2}y_1+\text{lin.comb.}\{z_1,x_2,y_2,z_2)({1\over\nu}x_2+\text{lin. comb.}\{y_2,z_2\}]=$$
$${1\over\nu^3}[y_1x_2]+\text{lin. comb.}\{[y_1][y_2],[x_2z_2],[y_2z_2]\}+\text{lin. comb.}\{y_1z_2,z_1x_2,z_1y_2,z_1z_2\}.$$
One can check that the latter linear combination cannot be a cocycle (if it is nontrivial).
This yields the necessary expression for the row of the matrix of $\eta^*$. 
\subhead Case $H^i, i>2$
\endsubhead
Note that in this case, all the cohomology classes are products of cohomology classes of lower degrees. Hence, the following ordering of vectors yields the triangular matrix of $\eta^*$: one orders elements inductively (i.e., elements of degrees 1 and 2 are already ordered) and then writes the indecomposable elements of maximal degree multiplied by elements of lower degrees ordered at the previous step. Obviously, this order yields the triangular form of matrix $\eta^*$.
\vskip6pt
Elements, satisfying condition $(**)$ of the previous section, are found by a straightforward computation. For example,
$$\eta^*([x_1z_1])=[(\nu x_1+...)({1\over\nu} z_1+...)]=[x_1z_1]+... .$$
Vectors $(**)$ are the only ones with this property, since the coefficient $\nu$ cannot equal $1$ (see the proof of Lemma 5).
\head 5. Addendum: methods of calculating   the minimal models of free algebras and bigraded models
\endhead
Let us start with the method of calculating the minimal model of any free differential graded algebra.
\vskip6pt
In this section we follow the notation of \cite{L}, and here  $(\Lambda V,d)$ denotes any free (not minimal) differential graded algebra (hence, it is not $\Lambda V$ from Section 3. Assume that $V=\oplus_{i\geq 1}$ and that $H^1(\Lambda V,d)=0$. Denote by $\Lambda^{++}V$ the ideal in $\Lambda V$ generated by decomposables. Consider the map $d': V^n\to V^{n+1}$ being a composition of $d: V^n\to V^{n+1}\oplus(\Lambda^{++}V)^{n+1}$ and a projection on $V^{n+1}$. Let $V'$ be a subspace in $\operatorname{Im}\,d'$ such that
$$V=\operatorname{Im}\,d'\oplus V'\oplus W,\,\, \operatorname{Ker}\,d'=\operatorname{Im}\,d'\oplus V'.$$
Then, $d'|_{W}$ is an isomorphism. Consider the image $W'=d(W)\subset \Lambda V$. 
\proclaim{Proposition}\cite{L, Section 8} 
\roster
\item There is an isomorphism of graded algebras
$$\Lambda V\cong \Lambda V'\otimes\Lambda(W'\oplus W)$$
where $C=\Lambda(W'\oplus W)$ is contractible, and $V'=H^*(V,d')$.
\item The minimal model of $\Lambda V$ is a factor algebra
$$(\Lambda V/\langle C^+\rangle,\bar d)$$
where $\langle C^+\rangle$ denotes an ideal in $\Lambda V$ generated by  elements of $C$ of positive degrees, and $\bar d$ is a differential induced by $d$.
\endroster
\endproclaim
 
\remark{Remark} Assumption $H^1(\Lambda V,d)=0$ is essential. One can show by examples (cf. \cite{BG, FLS, TO}) that without this condition the result does not hold. Note that $(\Lambda\,(Z_0^U\oplus Z_1^U\oplus Z_2^U\oplus...), D_U)$ constructed in Section 3, is connected and simply-connected.
\endremark
Now, we give a brief description of the procedure of constructing the bigraded model of any commutative graded algebra.
\vskip6pt 
Let $H$ be any commutative graded algebra regarded as a differential graded algebra with zero differential. There is a commutative differential graded algebra $(\Lambda\,Z,d)$ with the following properties:
\roster
\item $Z=\oplus_{n=0}^{\infty}Z_n$ (the lower grading)
\item $Z_{(n)}=Z_0\oplus ...\oplus Z_n$ yields a $d$-invariant filtration on $\Lambda\,Z$ and on $H^*(\Lambda\,Z)$,
\item $(\Lambda\,Z,d)$ is a minimal model of $(H,0)$.
\endroster
The method of constructing $(\Lambda\,Z,d)$ is given in \cite{HS, p.243} and goes as follows.
\subhead The space $Z_0$
\endsubhead
Set $Z_0=H^+/H^+\cdot H^+$, it is the space of indecomposables of $H$. Set $d=0$ in $Z_0$. Define $\rho: \Lambda\,Z_0\to H$ as a splitting of the projection $H^+\to Z_0$. Denote the kernel of the section by $K$.
\subhead The space $Z_1$
\endsubhead
Set $Z_1=K/(K\cdot \Lambda^+Z_0)$ with a shift downward by one of degrees:
$$Z_1^p=(K/K\cdot\Lambda^+Z_0)^{p+1}.$$
Then, extend $d$ to $Z_1$ by requiring that it be a linear map $Z_1\to K$ splitting the projection.
\subhead The spaces $Z_n$
\endsubhead Suppose $(\Lambda\,Z_{(n)},d)$ has been constructed for some $n\geq 1$ so that $d$ is homogeneous of lower degree $-1$. Define $Z_{n+1}$ by
$$Z_{n+1}^p=(H_n(\Lambda\,Z_{(n)},d)/H_n(\Lambda\,Z_{(n)},d)\cdot H_0^+(\Lambda\,Z_{(n)},d)^{p+1}.$$
\remark{Remark} In \cite{HS} one can find explicit examples of the inductive calculation above. Also, \cite{HS} contains examples which show that after perturbation of the differential, the resulting differential graded algebra may become not minimal, since linear terms may appear.
\vskip6pt
Note that one can think about $Z_0$ as the space of generators of $H$ and $Z_1$ as the "space" of relations between generators.
\endremark

\enddocument
\end